\documentclass[11pt]{article}

\usepackage{amsmath,fullpage,amssymb,graphicx}
\usepackage{fourier,eucal}
\usepackage[bbgreekl]{mathbbol}

\newcommand{\sia}[1]{{\mathsf{#1}}_i^\alpha}

\newcommand{\hil}[1]{\hat{#1}_i^\ell}

\newcommand{\il}[1]{{#1}_i^\ell}
\newcommand{\ik}[1]{{#1}_i^k}

\newcommand{\iko}[1]{{#1}_i^{k+1}}
\newcommand{\ia}[1]{{#1}_i^\alpha}
\newcommand{\bia}[1]{\bar{#1}_i^\alpha}
\newcommand{\bbia}[1]{\bar{\bar{#1}}_i^\alpha}
\newcommand{\bial}[1]{\bar{#1}_i^{\alpha,\ell}}
\newcommand{\bbial}[1]{\bar{\bar{#1}}_i^{\alpha,\ell}}

\newcommand{\tia}[1]{\tilde{#1}_i^{\alpha}}
\newcommand{\ttia}[1]{\tilde{\tilde{#1}}_i^{\alpha}}

\begin{document}

\begin{titlepage}
\title{\bf Finite Element Methods with Artificial Diffusion for Hamilton-Jacobi-Bellman Equations}
\author{Max Jensen\footnote{Durham University, Durham, United Kingdom, m.p.j.jensen@durham.ac.uk}, Iain Smears\footnote{Oxford University, Oxford, United Kingdom, iain.smears@maths.ox.ac.uk}}
\end{titlepage}

\maketitle

\begin{abstract}
In this short note we investigate the numerical performance of the method of artificial diffusion for second-order fully nonlinear Hamilton-Jacobi-Bellman equations. The method was proposed in (M.~Jensen and I.~Smears, arxiv:1111.5423); where a framework of finite element methods for Hamilton-Jacobi-Bellman equations was studied theoretically. The numerical examples in this note study how the artificial diffusion is activated in regions of degeneracy, the effect of a locally selected diffusion parameter on the observed numerical dissipation and the solution of second-order fully nonlinear equations on irregular geometries.
\end{abstract}

\section{Introduction} \label{sec:Introduction}

Hamilton-Jacobi-Bellman (HJB) equations, which describe the value function in the theory of optimal control, are fully nonlinear partial differential equations, which are of second-order if the underlying control problem is stochastic. One challenge arising in the numerical solution of these equations is the presence of spurious generalised solutions of the PDE which do not coincide with the value function. While these spurious solutions often possess the same regularity as the value function, they violate monotonicity properties exhibited by the value function. These properties lead to the concept of viscosity solutions. 

Regarding the numerical solution of HJB equations, we would like to highlight three approaches within the finite element methodology, which have been employed to ensure convergence to the value function. For a review of discrete Markov chain approximations before application of the dynamic programming principle we refer to \cite{KD01}. For the method of vanishing moments we point to \cite{FN11}. For the approach by Barles and Souganidis we cite the original source \cite{BS91} and the recent adaption \cite{JS11} by the authors to a finite element framework. A more comprehensive outline of the literature is in \cite{FGN11}.

\section{Numerical Method} \label{sec:Method}

Let $\Omega$ be a bounded Lipschitz domain in $\mathbb{R}^d$, $d \ge 2$. Let $A$ be a compact metric space and
\[
A \to C(\overline{\Omega}) \times C(\overline{\Omega}, \mathbb{R}^d) \times C(\overline{\Omega}) \times C(\overline{\Omega}),\; \alpha \mapsto (a^\alpha,  b^\alpha, c^\alpha, d^\alpha)
\]
be continuous. Consider the bounded linear operators of non-negative characteristic form 
\[
L^\alpha : \; H^2(\Omega) \cap H^1_0(\Omega) \to L^2(\Omega), \; w \mapsto - a^\alpha \, \Delta w + b^\alpha \cdot \nabla w + c^\alpha \, w
\]
where $\alpha$ belongs to $A$.
Then
\begin{align} \label{Lbounds}
\sup_{\alpha \in A} \| \, (a^\alpha,  b^\alpha, c^\alpha, d^\alpha) \,
\|_{C(\overline{\Omega}) \times C(\overline{\Omega}, \mathbb{R}^d) \times C(\overline{\Omega}) \times C(\overline{\Omega})} + 
\sup_{\alpha \in A} \| L^\alpha \|_{C^2(\overline{\Omega}) \to C(\overline{\Omega})} < \infty.
\end{align}
We assume that the final-time boundary data \( v_T \in C(\overline{\Omega}) \) is non-negative: \( v_T \geq 0 \) on \(\overline{\Omega}\). For smooth $w$ let $H w := \sup_\alpha (L^\alpha w - d^\alpha)$, where the supremum is applied pointwise. The HJB equation considered is
\begin{subequations}
\begin{alignat}{2}
- v_t + H v &=0 	&	&\qquad\text{in }(0,T)\times\Omega,\\
v&=g 		  	&	&\qquad\text{on }(0,T)\times\partial\Omega,\\
v&=v_T		  	&	&\qquad\text{on }\{T\}\times\overline{\Omega}.
\end{alignat}
\end{subequations}
Let $V_i$, $i=1,2,\dots$, be a sequence of piecewise linear shape-regular finite element spaces, whose underlying meshes are strictly acute. Let $y_i^\ell$, $\ell=1,\dots,\dim V_i$, denote the nodes of the mesh with associated hat functions $\phi_i^\ell$. Set $\hat{\phi}_i^\ell := \phi_i^\ell / \| \phi_i^\ell \|_{L^1(\Omega)}$. The mesh size is denoted $(\Delta x)_i$. The set of time steps is $S_i := \{ s_i^k : k = 0, \ldots, T / h_i \}$. Let the $\ell$th entry of $d_i w( s_i^k, \cdot)$ be $(d_i w( s_i^k, \cdot))_\ell = (w( s_i^{k+1}, y_i^\ell) - w(s_i^k, y_i^\ell)) / h_i$.

For each $\alpha$ and $i$ find an approximate splitting $L^\alpha \approx \ia E + \ia I$ into linear operators
\begin{align*}
\ia E &: \; \; w \mapsto - \bia a \, \Delta w + \bia b \cdot \nabla w + \bia c \, w,\\
\ia I &: \; \; w \mapsto - \bbia a \, \Delta w + \bbia b \cdot \nabla w + \bbia c \, w,
\end{align*}
of the form $a^\alpha = \tia a + \ttia a$, $b^\alpha = \bia b + \bbia b$, $c^\alpha = \bia c + \bbia c$ and $d^\alpha = d_i^\alpha$. To impose monotonicity, select the artificial diffusion parameters \( \bial \nu \) and \( \bbial \nu \) as in \cite{JS11} such that $\bia a(\il y) \geq \max \{ \tia a(\il y), \bial \nu \}$ and $\bbia a(\il y) \geq \max \{ \ttia a(\il y), \bbial \nu \}$.

Define, for $w \in H^1(\Omega)$, \(\ell \in \{ 1, \dots, N = \dim V_i^0 \} \),
\begin{subequations}\label{eq:discreteop}
\begin{align}
(\sia E w)_\ell & := \bia  a(\il y) \langle \nabla w, \nabla \hil \phi \rangle + \langle \bia  b \cdot \nabla w + \bia  c \, w, \hil \phi \rangle,\\
(\sia I w)_{\ell} & := \bbia a(\il y) \langle \nabla w, \nabla \hil \phi \rangle + \langle \bbia b \cdot \nabla w + \bbia c w, \hil \phi \rangle,\\
(\sia C)_\ell & := \langle \ia d, \hil \phi \rangle.
\end{align}
\end{subequations}
Obtain the numerical solution $v_i(T, \cdot) \in V_i$ by interpolation of $v_T$. Then $v_i(\ik s, \cdot) \in V_i$ at time $\ik s$ is defined, inductively, by interpolating the boundary data and by
\begin{align} \label{num}
- d_i v_i(\ik s,\cdot) + \sup_\alpha \bigl( \sia E v_i(\iko s,\cdot) + \sia I v_i(\ik s,\cdot) - \sia C \bigr) = 0,
\end{align}
where the supremum is understood to be applied component-wise to the collection of vectors
\[
\left\{\sia E v_i(\iko s,\cdot) + \sia I v_i(\ik s,\cdot) - \sia C \,:\,\alpha\in A\right\}.
\]

\section{Selection of the Artificial Diffusion Parameter} \label{sec:ad}

That the diffusion coefficients $\bia{a}$ and $\bbia{a}$ in \eqref{eq:discreteop} are placed outside of the scalar product originates from the non-divergence form of the linear operators $L^\alpha$ of HJB operators \cite{JS11}. This structure makes it straightforward to implement {\em local} artificial diffusion parameters \( \bial \nu \) and \( \bbial \nu \), that is to implement a dependence on the nodal position $\il y$. Indeed a change of  \( \bial \nu \) corresponds in the assembly of $\sia E$ to a scalar multiplication of the $\ell$th row of the stiffness matrix, see \eqref{eq:discreteop}. 

Changes in the optimal artificial diffusion coefficients arise from variations of the local mesh quality and the local mesh P\'eclet number. The parameters \( \bial \nu \) and \( \bbial \nu \) may be chosen by studying (35) in \cite{JS11} or by examining the assembled matrices of the unstabilised operators. This is in particular simple for the $\sia E$, since only the signs of off-diagonal entries need be corrected to impose a local monotonicity property, thus leading to an algorithm which can be executed row-by-row.

Notice that in contrast, for problems in divergence form, the stabilised diffusion coefficients $\tia{a}$ and $\ttia{a}$ need to be determined in the whole domain and not just at nodal positions. 

\section{An Exact Solution and Convergence Rates} \label{sec:Rates}

\begin{figure}[t]
\begin{center}
 \includegraphics[width=8cm]{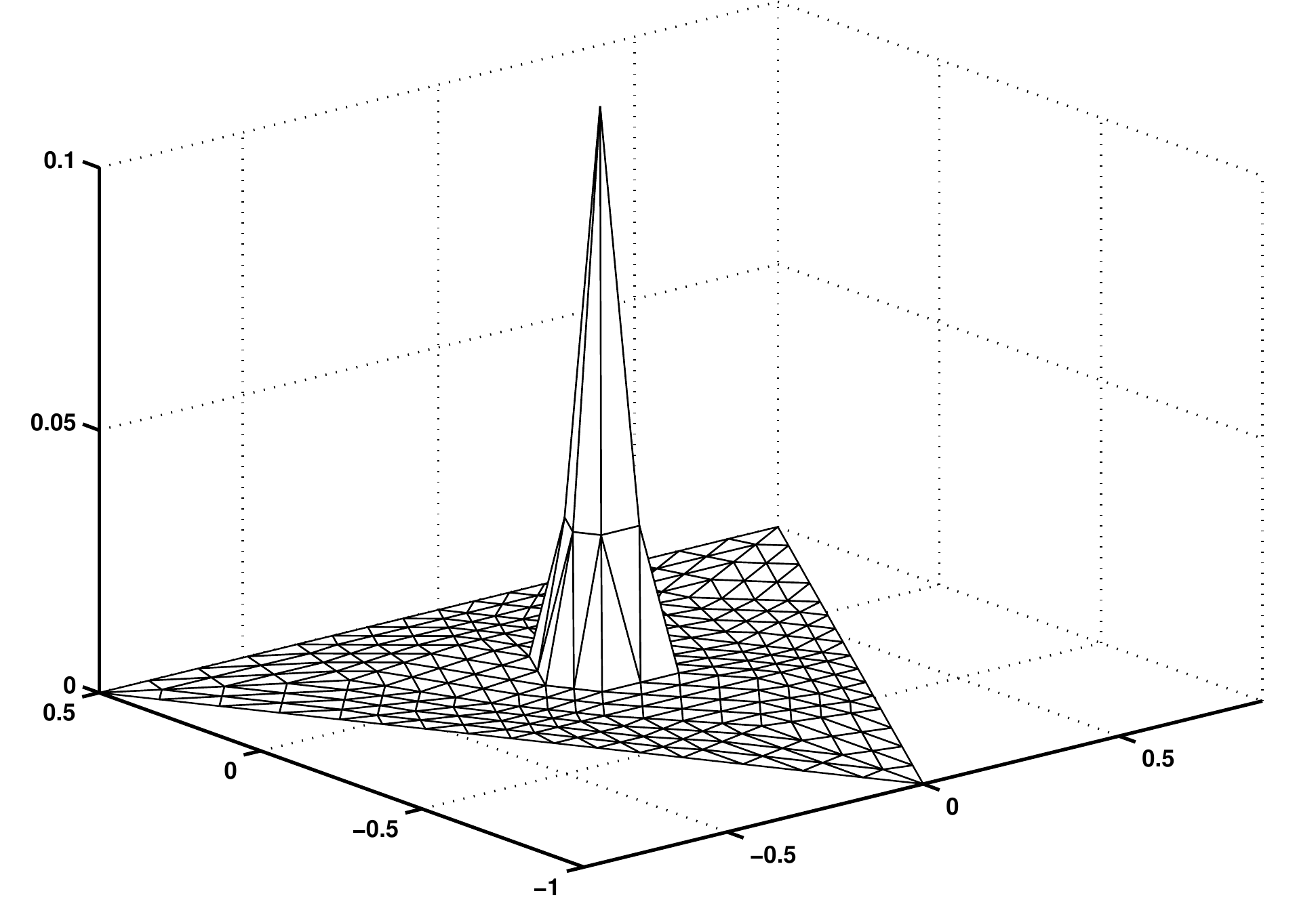}
\caption{The plot shows a locally-adapted choice of artificial diffusion over the triangular domain. The peak in artificial diffusion at the centre of the domain corresponds to the degeneracy of the differential operator at the origin.}
\label{fig:art_diff_triangle}
\end{center}
\end{figure}

We consider a triangular spatial domain $\Omega$ with the vertices $(0,-1)$, $(\sqrt{3}/2,1/2)$ and $(-\sqrt{3}/2,1/2)$ and the HJB equation
\begin{align} \label{eq:HJBtri}
 -v_t - \frac{1}{2} \sqrt{\frac{x^2+y^2}{T-t+1}} \Delta v + \frac{1}{2} \frac{1}{\sqrt{T-t+1}} \left| \nabla v \right| = - \frac{1}{2} \frac{\sqrt{x^2+y^2}}{(T-t+1)^{3/2}}.
\end{align}
To see that equation \eqref{eq:HJBtri} is a HJB equation, note that the Euclidean norm of the gradient satisfies
\begin{align} \label{eq:eik}
\left| \nabla v \right| = \sup \{ \beta \cdot \nabla v\, :\, \beta \in \mathbb{R}^2 \text{ with } \left|\beta\right| = 1 \}.
\end{align}
A calculation verifies that the function
\[
 v(x,y,t) = \exp\left( - \sqrt{\frac{x^2+y^2}{T-t+1}} \right) + \sqrt{\frac{x^2+y^2}{T-t+1}}
\]
is an exact solution of \eqref{eq:HJBtri}, where boundary and final-time conditions are determined by interpolation of the exact solution $v$. The equation is to be solved on the time interval $(0,1)$.

\begin{figure}[t]
\begin{center}
\includegraphics[width=11cm]{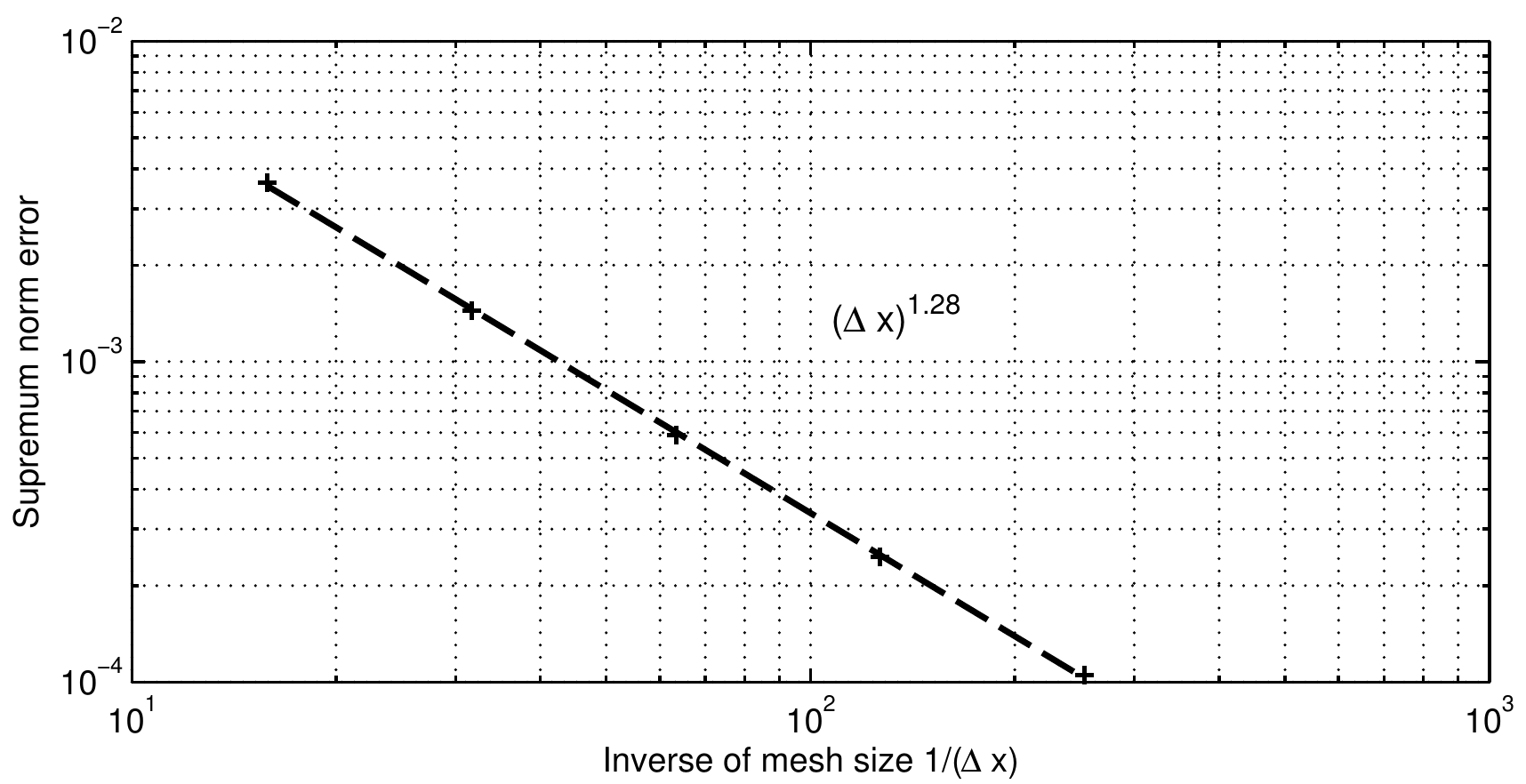}
\caption{Uniform error at the initial time $t=0$ with constant ratio $h_i / (\Delta x)_i$.}
 \label{fig:stab}
\end{center}
\end{figure}

We consider a splitting which discretises the advection term explicitly with the minimum amount of diffusion needed for monotonicity. The remaining (linear) diffusion is incorporated in the implicit term, observing that this leads to a time-step restriction $h_i \lesssim  (\Delta x)_i$. Figure~\ref{fig:stab} shows the supremum norm error with a constant ratio $h_i / (\Delta x)_i$, showing the uniform stability of the method in this setting.

\begin{figure}[t]
\begin{center}
\includegraphics[width=10.75cm]{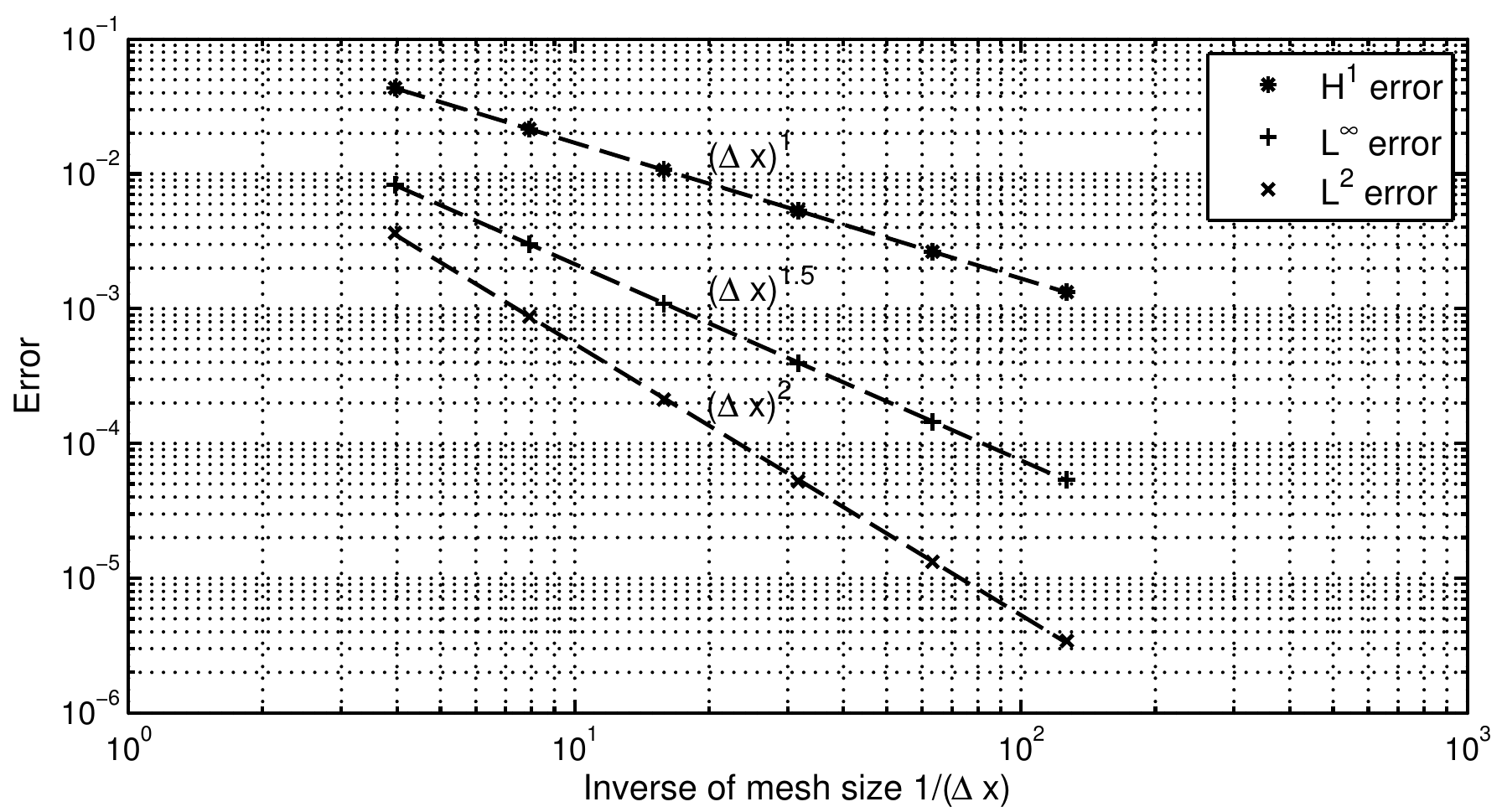}
\caption{Convergence rates of the numerical solution at the initial time $t=0$ in the $L^2$-, $L^\infty$- and $H^1$-norms with constant ratio $h_i / (\Delta x)_i^2$.}
 \label{fig:rates}
\end{center}
\end{figure}

Figure~\ref{fig:art_diff_triangle} illustrates a choice of artificial diffusion that is locally adapted to the P\'{e}clet number of a coarse mesh. It is seen that diffusion is only artificially introduced in a neighbourhood of the origin, which is where the operator becomes degenerate. Figure~\ref{fig:rates} illustrates the rates of convergence of the numerical solution $v_i(0,\cdot)$ in the $L^2$-, $L^\infty$- and $H^1$-norms to the exact solution at the initial time, now using the constant ratio $h_i / (\Delta x)_i^2$.
  
\section{Eikonal Equation} \label{sec:Eikonal}

\begin{figure}[t]
\begin{center}
\includegraphics[width=11.8cm]{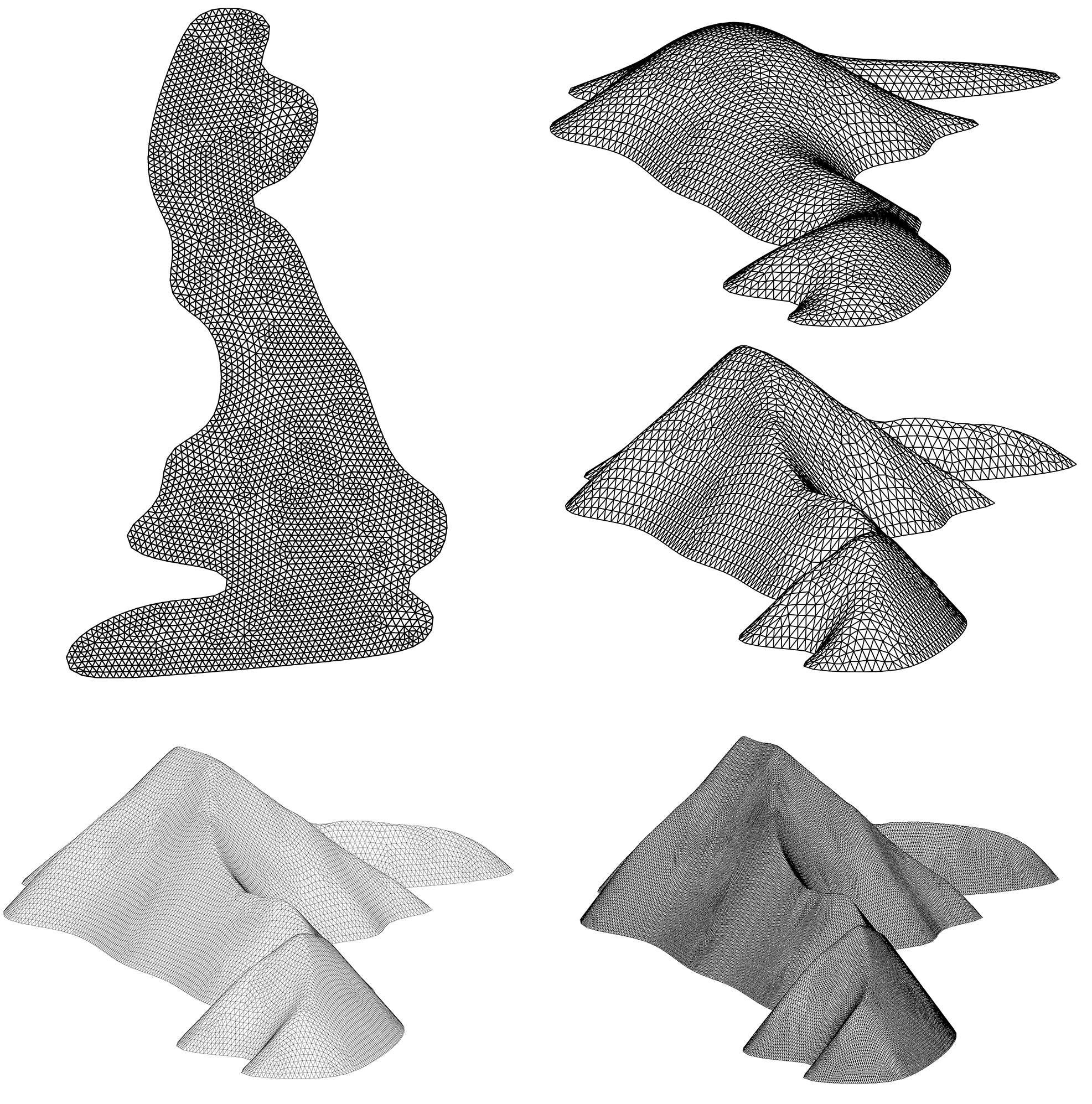}
\caption{The top left plot shows the original mesh for the domain. The top right plot shows the numerical solution with globally-chosen diffusion and below is a plot with locally-chosen artificial diffusion. The bottom left plot is the numerical solution from one uniform refinement of the mesh and the bottom right plot shows the solutions after two refinements.}
\label{fig:hj_surf}
\end{center}
\end{figure}

The steady-state limit of the time-dependent eikonal equation,
\[
 -v_t + \left| \nabla v \right| = 1,
\]
equipped with homogeneous boundary and final-time conditions, measures the distance to the boundary. Due to \eqref{eq:eik} the eikonal equation belongs to the Hamilton-Jacobi-Bellman family. We consider the equation on a domain (Figure~\ref{fig:hj_surf}, top left) whose irregular shape leads to complicated curves on which $v$ is not differentiable. The height of the domain is equal to one.

Figure \ref{fig:hj_surf} compares locally-adapted choices of the artificial diffusion parameter with global choices, and illustrates the effect of mesh refinement on the quality of the solution. To compare the quality of the numerical solutions we compare their $L^\infty$-norms---recalling that excessive numerical dissipation leads to a smearing out of extrema. The coarse grid solution, with $2858$ internal nodes and with a global diffusion parameter of $0.05$, only reaches a height of $9.56\times10^{-2}$. In contrast, a computation on the same mesh with a local diffusion parameter, with mean value $0.01$ and standard deviation $0.005$ but maximal value $0.05$, leads to a height of $0.138$. With one (two) steps of uniform refinement the number of elements increases by a factor of $4$ (of $16$) and the artificial diffusion decreases to an average value of $0.004$ (of $0.002$), giving the improved value of $0.147$ (of $0.151$) for the $L^\infty$-norm in the lower left (right) part of Figure \ref{fig:hj_surf}. The $L^\infty$-norm of a reference solution on a very fine mesh is $0.153$. 

\section{A Second-Order Fully Nonlinear Equation} \label{sec:Fully}

\begin{figure}[t] 
\begin{center}
\includegraphics[height=6cm]{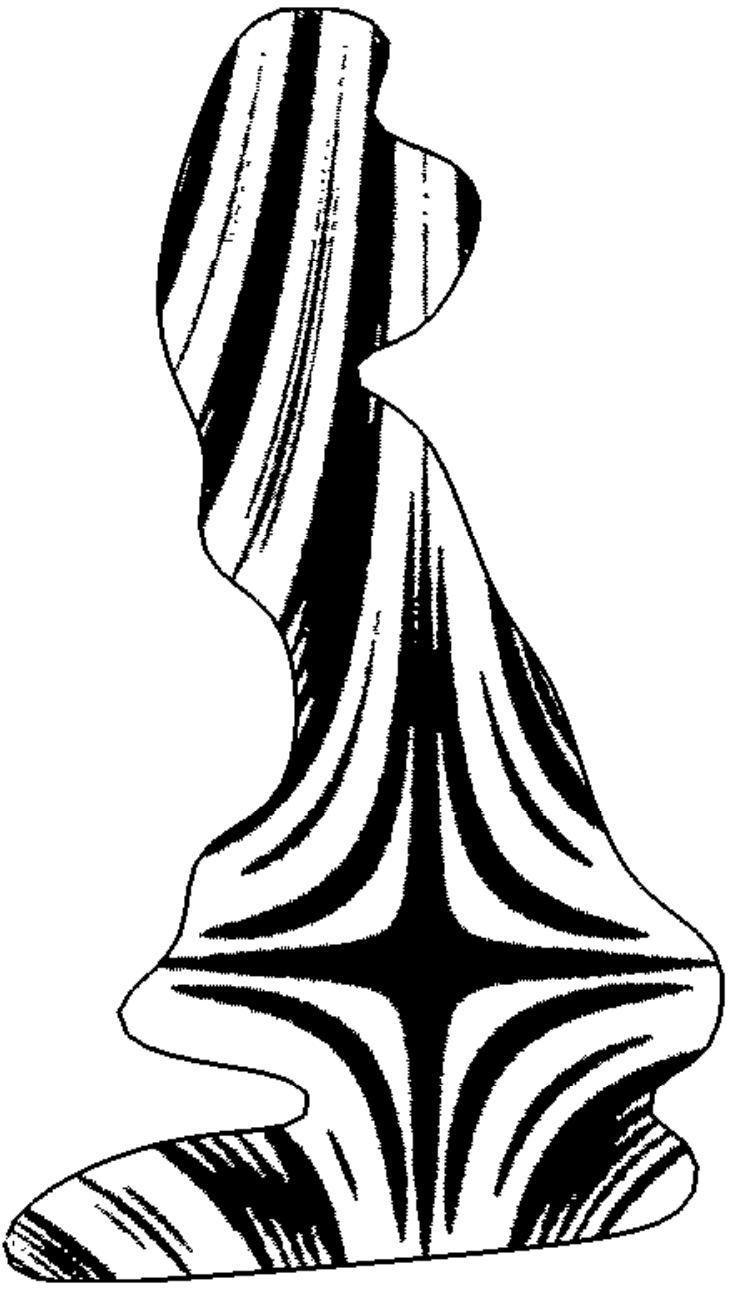} \hspace{10mm}
\includegraphics[height=6cm]{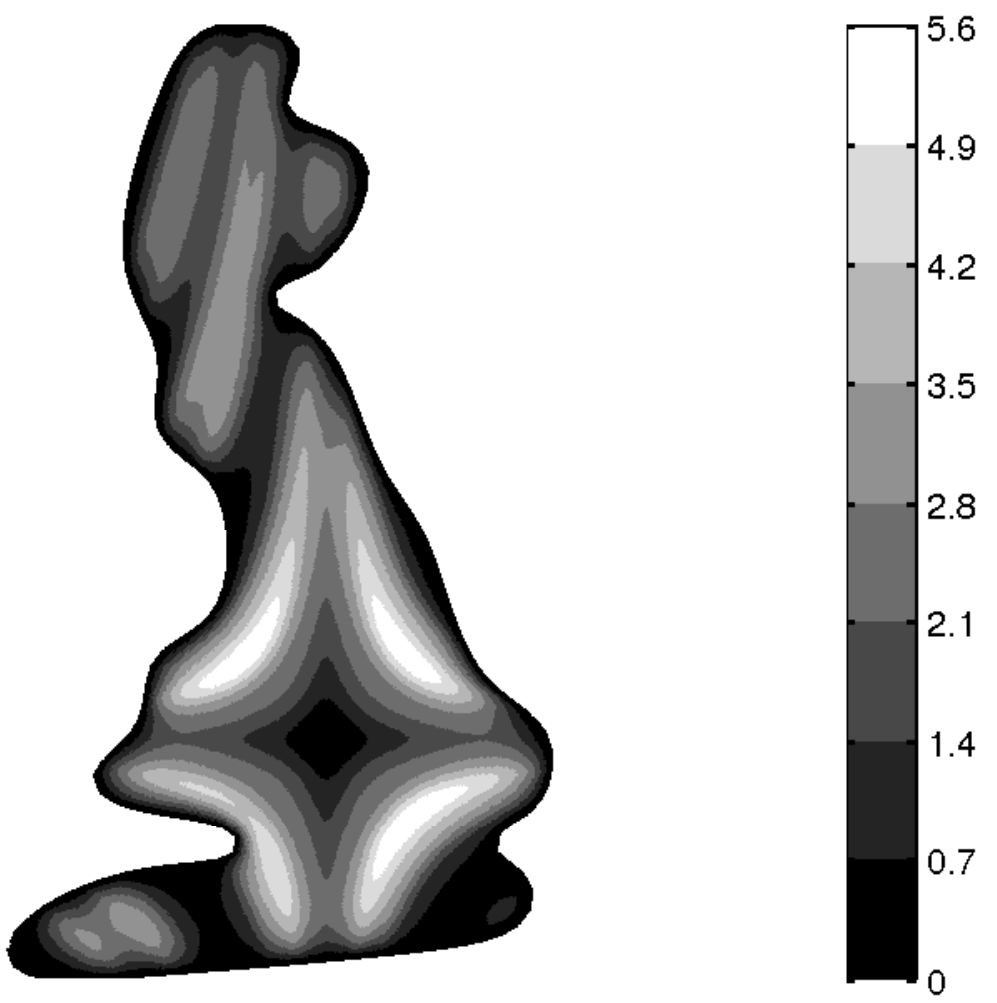}
\caption{The left plot shows the contours of the control $\alpha$ maximising the nonlinear second order term of the equation, the right plot shows the value function $v$. In the left plot, the black regions correspond to $\alpha = \alpha_0$, whilst the white regions correspond to $\alpha = \alpha_1$.}
\label{fig:fully}
\end{center}
\end{figure}

The final example, on the same domain as in the previous section, is concerned with the fully nonlinear case of second-order. We examine the equation
\[
 -v_t + \sup_{\alpha \in [\alpha_0,\alpha_1]}\{ -\alpha \Delta v \} + \left|\nabla v\right| = -v_t + \sup_{(\alpha,\beta) \in [\alpha_0,\alpha_1]\times \partial B(0,1)}\{ -\alpha \Delta v  + \beta\cdot\nabla v \} = f,
\]
where $B(0,1)$ is the unit ball in $\mathbb{R}^2$, $T = 0.009$, $\alpha_0 = 0.045$, $\alpha_1 = 0.09$ and
\begin{align*}
 f(x,y) &=529 \biggl( \sin\bigl( g(x,y) \bigr) + \frac{1}{2} \sin\bigl( 2\, g(x,y) \bigr) + \frac{4}{10} \sin\bigl(8\, g(x,y)\bigr) \biggr)^2,\\
 g(x,y) &= \pi^2 (x-0.63)(y-0.26)/0.07.
\end{align*}
The boundary and final-time conditions are homogeneous. As before, the advection term is discretised explicitly, with the locally minimal diffusion needed for monotonicity. The possibly remaining (fully nonlinear) diffusion is placed in the implicit term.

In this example the control $\alpha$ of the second-order term is maximised independently of the first-order control $\beta$; in the sense that the optimal $\alpha$ may be determined without knowledge of $\beta$. Furthermore, as $\Delta v$ takes locally either a positive or negative value only the controls $\alpha_0$ and $\alpha_1$ are ever active in the HJB equation. This is an example of the Bang-bang principle. It is also reflected in the left plot in Figure \ref{fig:fully}, where the value of $\alpha$ is plotted. Black colouring signifies $\alpha_0$ maximises the operator, whereas white colouring corresponds to $\alpha_1$. Observe that no intermediate grey values can appear. The plot of the control $\alpha$ mimics some of the features of the value function, which is plotted in the same figure on the right---in part because the Laplacian contains information about the curvature of the solution.

At each time-step of the method, a semi-smooth Newton method \cite{JS11} was used to solve the nonlinear discrete equation in \eqref{num}, where each iteration of the algorithm involves solving a linear system. To study the performance of the algorithm, the HJB equation was solved on a sequence of successively refined meshes, with a constant set of tolerances. The linear systems were solved to a tolerance of $10^{-10}$ by GMRES. The stopping criterion for Newton's method was a relative residual tolerance of $5 \times 10^{-8}$ in the maximum norm and a convergence of the iterations requirement of $5 \times 10^{-9}$ in the maximum norm. The sizes of the linear systems for the sequence of meshes were $674$, $2858$, $11759$, $47693$ and $192089$. The respective average number of Newton iterations for a single time step were $3$, $3.67$, $4.04$, $4.22$ and $4.86$, with respective standard deviations $0$, $0.48$, $0.19$, $0.42$ and $0.36$. This demonstrates a weak dependence of the number of iterations needed for an individual time step on the system size, thus showing that the total number of linear systems to be solved for a complete computation depends principally on the number of time-steps.


\bibliographystyle{spmpsci}

\end{document}